\documentclass[12pt,fleqn,a4paper]{article}
\usepackage{dsfont}
\usepackage{amssymb}
\usepackage{amsfonts}
\usepackage{amsmath}
\usepackage[english]{babel}     
\usepackage[pdftex,pdfstartview=FitH]{hyperref}
\def\hxi{\hat{\xi}}
\def\cost{{\rm cost}}
\def\comp{{\rm comp}}
\def\e{{\rm e}}
\def\P{{\rm P}}
\def\N{{\bf N}}

\def\rand{{\rm rand}}
\def\quant{{\rm quant}}
\def\det{{\rm det}}
\def\vsn{\vskip 1pc \noindent}
\def\ve{\varepsilon}
\def\res{{\rm res}}
\def\abs{{\rm abs}}

\def\fr{F_d^{r,\rho}}
\def\frr{F^{r,\rho}}
\def\fdrr{F_d^{r,\rho}}
\def\tfrr{\tilde{F}^{r,\rho}}
\def\R{{\bf R}}

\def\be{\begin{equation}}
\def\ee{\end{equation}}

\newtheorem{thm}{Theorem}

\renewcommand{\v}{\phi}
\newenvironment{pf}{\noindent \textit{Proof.}}{  \ \hspace{\stretch{1}} $\Box$\\ }
\newcommand{\intl}{\int\limits}

\newenvironment{pf*}[1]{\noindent\textit{#1} }{}

\begin{document}

\thispagestyle{empty}

\begin{center}
{\Large{\bf On the quantum complexity of approximating median of continuous distribution}}\footnotemark[1] \\
\vspace{1cm}
{\large Maciej~Go\'cwin\footnotemark[2]}
\end{center}

\footnotetext[1]{ ~\noindent This research was partly supported by the 
Ministry of Science and Higher Education
}

\footnotetext[2]{ 
\begin{minipage}[t]{16cm} 
 \noindent
{\it Faculty of Applied Mathematics, AGH University of Science 
and Technology,\\ 
\noindent  Al. Mickiewicza 30, paw. A3/A4, III p., 
pok. 301,\\
 30-059 Krakow, Poland }
\newline
E-mail: $\;\;$ gocwin@agh.edu.pl \\

\end{minipage} }

\vspace{2cm}
\begin{center}
{\Large{\bf Abstract}}
\end{center}
\vsn
We consider approximating of the median of absolutely continuous distribution given by a probability density function $f$.  We assume that $f$ has $r$ continuous derivatives, with  derivative
of order $r$ being H\"older continuous with the exponent $\rho$.
We study the $\ve$-complexity of this problem in the quantum setting.
We show that the $\ve$-complexity up to logarithmic factor is of order  $\ve^{-1/((r+\rho+1))}$. 
\\ \\
%{\bf Keywords:} 
%nonlinear equations, deterministic algorithms, randomized algorithms,
%quantum algorithms, optimality, complexity 

\newpage
\section{Introduction}\label{sec:introduction} The quantum complexity of discrete problems is well studied. A speed-up of quantum algorithms over deterministic and randomized algorithms is shown for many
problems, starting from the factorization algorithm of Shor \cite{Shor}, followed by
database search algorithm of Grover \cite{Grover}. Other discrete problems, such
as discrete summation, computation of the mean, median and
$k$th-smallest element were also studied, see e.g.
\cite{Brassard},\cite{Brassard2}\cite{Durr},\cite{Grover2},\cite{Nayak}.
 \\
There is also a progress in studying the quantum complexity of continuous analogues of these problems. The first paper dealing with
the quantum complexity of a continuous problem was the work of Novak
\cite{Novak1}, where the integration of a function from the
H\"older class is considered.  Integration in other
function spaces (Lebesgue and Sobolev classes) space was also
investigated and a quadratic speed-up over the randomized setting
was shown (see \cite{Heinrich},\cite{Heinrich3}). The problem
of function approximation on a quantum computer was studied by
Heinrich \cite{Heinrich1},\cite{Heinrich2}. The problem of maximization of function form the H\"older class was investigated in~\cite{Gocwin} and the problem of finding the root of the function was studied in~\cite{GocwinKacewicz}. Also path
integration \cite{Wozniakowski} and differential equations
\cite{Kacewicz2005} on a quantum computer were
investigated, and a speed-up was established.
\\
In this paper we consider the problem of approximating the median of continuous distribution. This problem is a continuous equivalent of the problem of finding the median of the discrete sequence which is investigated in~\cite{Nayak}. We assume that the distribution is given by density function which belongs to a
classical H\"older class. We present almost matching upper and lower complexity bounds in the
quantum setting. We show that quantum computations yield a
speed-up compared to deterministic and randomized
algorithms over entire range of class parameters.
\\
We extend also the problem to the problem of computing of the vector of quantiles of the continuous distribution and show matching upper and lower complexity bounds for this problem.
\\ 
In Section \ref{sec:formulation} the problem of computing the mean is formulated and  necessary definitions are presented. The complexity bounds for this problem are presented in Section~\ref{sec:median}. Section~\ref{sec:quantiles} contains the results for the problem of computing of the vector of quantiles.
%\section{Computing of the median}\label{sec:median}
\section{Problem formulation and basic definitions}\label{sec:formulation}
Let $f:[0,1]\rightarrow \R$ be a real-valued, nonnegative function which integrates to 1.  We are interested in approximation the median of the probabilistic distribution with probabilistic density function $f$  with precision $\ve>0$
in the sense of the absolute or the residual error criterion. Let $\xi:\ \intl_0^{\xi}f(t)dt=1/2$ be the median. We wish to compute a point $\hxi\in [0,1]$ for which $|\hxi-\xi|\leq \ve$  (the absolute criterion), or 
$\left|\intl_{0}^{\hxi}f(t)dt-\frac{1}{2}\right|\leq\ve$ (the residual criterion). We consider a H\"older class of functions $f$ given by 
\begin{eqnarray*}
\frr&=&\left\{f:[0,1] \rightarrow \R \,\, |\,\, f\in C^r([0,1]), |f^{(i)}(x)|\leq D,\,\, i=0,1,
\ldots r,\right.\\  & & \left.|f^{(r)}(x)-f^{(r)}(y)|\leq H |x-y|^{\rho} \;\mbox{ for } x,y\in[0,1] \right\},
\end{eqnarray*}
where $r\geq 0$, $0<\rho\leq 1$ and  $D,H$ are positive constants.  
Let us define the class of density function from the H\"older class $\frr$ by $\frr(1)$, that is
\[
\frr(1)=\left\{f\in\frr\,\,|\,\, f(x)\geq 0\ \forall x\in[0,1],\ \intl_0^{1}f(t)dt=1\right\}.
\]
For the absolute criterion we will also need a technical assumption that the function $f$ is separated form zero. Let
\[
\tfrr(1)=\left\{f\in\frr(1)\,\,|\,\, f(x)\geq \gamma\ \forall x\in[0,1]\right\},
\]
where $\gamma$ is a positive constant.\\ 
We shall analyze the problem of approximation the median of the distribution with density function $f\in\frr(1)$ or $f\in\tfrr(1)$ in the quantum setting and for the comparison in the deterministic worst-case and in the randomized  settings.\\
To find the approximation $t^*$ we need some information about the function $f$. 
In the deterministic setting, by
information we mean an operator $N$ that assigns to the function $f$ a vector $N(f)$ 
of length at most $n$, whose components are
the values of $f$ or its derivatives of order at most $r$ at some (adaptively chosen) deterministic points form $[0,1]$.
Let $(\Omega,\Sigma, \P)$ be a probability space.
In the randomized setting, by information we mean a family of operators $N=\{N^\omega\}_{\omega\in \Omega}$.
Here, $N^\omega(f)$ is defined by the values of $f$ or its derivatives computed at randomly chosen points
from $[0,1]$. We assume that the successive evaluation points and the derivative numbers are chosen adaptively as functions of information values computed so far. 
The decision whether to compute a successive piece of information is also taken based on the computed  values,
according to some termination criterion. 
We assume that the length $n^\omega(f)$ 
of the information vector  is a measurable function and it satisfies 
${\rm E}(n^{\omega}(f))\leq n$ for all $f\in \fr$.\\
In the quantum setting, information is gathered by applying $n$ times a (standard) quantum query operator
 about the function 
$f$. For a detailed description of the quantum model of computation, quantum algorithms, and the quantum query operator, the reader is referred to \cite{Heinrich, Nielsen}. Roughly speaking, a quantum query is a unitary operator 
defined by $f$ whose application plays a role of evaluating the value of a function or its derivative
in the standard deterministic setting. \\
The approximate solution $\hxi$ is obtained in the deterministic setting by an algorithm $\v$ that computes  $\hxi$ based on the information values $N(f)$, $\hxi=\v (N(f))$. 
In the randomized setting, the algorithm is a family $\v=\{\v^{\omega}\}_{\omega\in \Omega}$,
and $\hxi=\hxi^{ (\omega)}=\v^{\omega}(N^{\omega}(f))$.
In the quantum setting, 
the algorithm is defined by an application of a number of sequences of unitary operations, $n$ of them being quantum queries,
each of the sequence followed by a quantum measurement. Further possible operations performed on a classical computer 
lead to $\hxi=\hxi^{(\omega)}$, see \cite{Heinrich,Nielsen}. 
In the randomized and quantum settings, the point $\hxi=\hxi^{(\omega)}$ 
is a random number.\\
We shall consider the absolute and residual error criteria. 
The superscript '$\omega$' in the definitions below in the deterministic setting is obviously irrelevant,
 and it will then be omitted.\\
%Let $S(f)=\{x\in [0,1]^d:\,\, f(x)=0\}$.
We define the local absolute error of an algorithm $\v$ by
\[
\e_{\rm abs}^{\omega}(f,\v^\omega)=\left|\hxi-\xi\right|.
\]
The residual error is given by
\[
\e_{\res}^{\omega}(f,\v^\omega)=\left|\intl_0^{\hxi}f(t)dt-1/2\right|.
\]
 In the randomized and quantum settings, we 
assume that $\e_{\rm abs}^{\omega}(f,\v^\omega)$ and $\e_{\res}^{\omega}(f,\v^\omega)$ are random variables.\\
Let the subscript '${\rm crit}$' stand for '${\rm abs}$' or '${\rm res}$'.
For $F$ being $\frr(1)$ or $\tfrr(1)$ the global error in the class $F$ in the deterministic, randomized and quantum settings
 is defined respectively by:\\
-- in the deterministic setting
\[
\e^{\det}_{\rm crit}(F,\v)=\sup_{f\in F}\e_{\rm crit}(f,\v),
\]
-- in the randomized setting
\begin{equation}\label{e_rand}
\e^{\rand}_{\rm crit}(F,\v)=\sup_{f\in F}\left(\intl_{\Omega}(\e_{\rm crit}^{\omega}
(f,\v^\omega))^2\,d\P(\omega)\right)^{1/2},
\end{equation}
-- in the quantum setting 
\[
\e^{\quant}_{\rm crit}(F,\v)=\sup_{f\in F}\; \inf\{\alpha:
\P(\,\e_{\rm crit}^{\omega}(f,\v^\omega)>\alpha \,)\leq 1/4\}.
\]
%\newline
Hence, in all cases the worst case with respect to functions $f$ from $F$ is of interest. 
In the randomized setting the error is defined as usual by the expected value, while
in the quantum setting by the demand to achieve the success probability of $\v$ at least 3/4. \\
The cost of an algorithm $\v$ in the deterministic and randomized settings is defined as an (expected)
 number of 
computed information values, that is a number of evaluations of function $f$ or its derivatives 
which an algorithm is based on. In the quantum setting, the (query) cost is meant as the number of 
applications of the quantum query operator. We denote the cost of an algorithm $\v$ 
in the class $F$ by $\cost(F,\v)$ with suitable superscripts: '$\det$', '$\rand$' or '$\quant$'
in the proper setting. 
\newline
For $\ve>0$, the $\ve$-complexity in the class $F$ is defined 
as a minimal cost of an algorithm that approximate the median of a distribution  with a density function from the class $F$ with the precision 
at most $\ve$,
\[
 \comp_{\ve,{\rm crit}}^{\#}(F)=\inf_{\v}\{\, \cost^{\#}(F,\v):\;\; \e^{\#}_{\rm crit}
(F,\v)\leq \ve \, \},
\]
where $ \#\in\{\det,\rand,\quant\}$. In the quantum setting, the complexity defined above is the 
quantum query 
complexity (the number of qubits can be arbitrary). 
It is also possible to consider, for instance, the qubit complexity, defined as the minimal
number of qubits sufficient to achieve the error at most $\ve$ (here the number of queries 
can be arbitrary), or other mixed type complexities.\\
Our aim is to show possibly tight upper and lower bounds on $\comp_{\ve}^{\#}(F)$.
Firstly, we recall in the next section a result on the problem of summation, integration and  solving initial value problems.

\section{Useful results}\label{sec:known}
We will recall in this section the useful known results, that will help us to prove the complexity bounds on the problem of approximating the mean and quantiles of continuous distribution.\\
We start with the problem of computing the mean of a sequence. We will use these results to prove the lower bounds. Let $g:\{0,1,\ldots,N-1\}\rightarrow [0,1]$ be a discrete function. We are interested in approximation of the mean of the sequence $(g(0),g(1),\ldots,g(N-1))$, that is the number $\displaystyle \frac{1}{N}\sum_{i=0}^{N-1}g(i)$ with precision $\ve>0$. Consider the deterministic, randomized and quantum settings. The definition of the error,  the cost and the complexity are similar as in Section~\ref{sec:formulation}. Denote the $\ve$-complexity of this problem by $\comp_{\ve}(N)$ with a suitable superscript indicating the setting.  From the known complexity bounds for this problem (see \cite{Novak,Mathe} for the deterministic and the randomized setting, \cite{Grover2,Nayak} for the quantum setting) we have that 
\begin{itemize}
\item in the deterministic setting   \begin{equation}\label{lower_detst}\comp_{\ve}^{\det}(N)=\Theta (N(1-2\ve));\end{equation} 
\item in the randomized setting \begin{equation}\label{lower_rand}\comp_{\ve}^{\rand}(N)=\Theta(\min\{N,\ve^{-2}\});\end{equation}
\item in the quantum setting \begin{equation}\label{lower_quant}\comp_{\ve}^{\quant}(N)=\Theta(\min\{N,\ve^{-1}\}).\end{equation}
\end{itemize}

We will need also the upper bounds on the complexity of approximating the integral, that we will use to show the upper bounds. Let $g:[0,1]\rightarrow \R$ be a function from H\"older class $\frr$. Our aim is to approximate $\int_0^1g(x)dx$ with precision $\ve$. Denote the complexity of this problem by $\comp_{\ve}(Int,\frr)$ (with suitable superscript). It is known (see e.g.~\cite{Novak} for the deterministic and randomized settings and \cite{Novak1} for the quantum setting), that the complexity of this problem is bounded by 
\begin{itemize}
\item in the deterministic setting   \begin{equation}\label{int_det}\comp_{\ve}^{\det}(Int,\frr)=O(\ve^{-1/{(r+\rho)}});\end{equation} 
\item in the randomized setting \begin{equation}\label{int_rand}\comp_{\ve}^{\rand}(Int,\frr)=O(\ve^{-1/{(r+\rho+1/2)}});\end{equation}
\item in the quantum setting \begin{equation}\label{int_quant}\comp_{\ve}^{\quant}(Int,\frr)=O(\ve^{-1/{(r+\rho+1)}}).\end{equation}
\end{itemize}
In the randomized setting, besides of the average error, we will also need to know the distribution of the error. Let $\phi$ be the optimal randomized algorithm. Then, there exists positive constant $C$ such that for any $g\in\frr$, $\ve>0$ and $\delta\in(0,1/2)$ we have
\[
\P\left(\left|\int_0^1g(x)dx-\phi(g)\right|>\ve\right)\leq\delta
\]
with cost 
\begin{equation}
\cost^{\rand}(\phi)\leq C \ve^{-1/{(r+\rho+1/2)}}\log(1/\delta).
\end{equation}

We will also need some results on the initial-value problems in ordinary differential equations.  The problems have the form
\begin{equation}
\label{problem_initial}
\left\{ \begin{array}{l}
z'(x)=f(z(x)),\qquad x\in[a,b]\\
z(a)=\eta
\end{array}\right.,
\end{equation}
where $a<b$, $f:\R^d\rightarrow \R^d$, $z:[a,b]\rightarrow \R^d$, and
$f$ belongs to the H\"older class
{\setlength\arraycolsep{0pt}\begin{eqnarray*}
\fdrr=\{&&f:\R^d\rightarrow \R^d|\ f\in C^r(\R^d), |D^{(i)}f(x)|\leq D,\ i=0,1,\ldots r,\\ 
&&\|D^{(r)}f(x)-D^{(r)}f(y)\|\leq H\|x-y\|^{\rho},\ \ x,y\in\R^d\},
\end{eqnarray*}}where $r\geq 0$, $\rho\in(0,1]$, $D$ and $H$ are positive constants. Denote the $\ve$-complexity of this problem by $\comp_{\ve}(IVP,\fdrr)$. There are known complexity bounds of this problem (see \cite{Kacewicz87} for the deterministic setting, \cite{Heinrich_Milla,Daun} for the randomized setting and \cite{Kacewicz2005} for the quantum setting):
\begin{itemize}
\item in the deterministic setting   \begin{equation}\label{ivp_det}\comp_{\ve}^{\det}(IVP,\frr)=O(\ve^{-1/{(r+\rho)}});\end{equation} 
\item in the randomized setting \begin{equation}\label{ivp_rand}\comp_{\ve}^{\rand}(IVP,\frr)=O(\ve^{-1/{(r+\rho+1/2)}});\end{equation}
\item in the quantum setting \begin{equation}\label{ivp_quant}\comp_{\ve}^{\quant}(IVP,\frr)=O(\ve^{-1/{(r+\rho+1-\gamma)}}),\end{equation}
\end{itemize}
for arbitrary $\gamma\in(0,1)$ (constant in big-O notation in the quantum setting may depend on $\gamma$).

\section{Complexity bounds}\label{sec:median}
\subsection{Upper complexity bounds}\label{sec:upper}
First we will present a theorem that states the upper complexity bounds on the approximation of the median of absolutely continuous distribution in the deterministic, randomized and quantum setting for both: the absolute and the residual criterion.
\begin{thm}
There exist positive constants $C^{\det}_{\abs}$, $C^{\rand}_{\abs}$ and $C^{\quant}_{\abs}$ which depend on $r$, $\rho$, $D$ and $\gamma$, and $C^{\det}_{\res}$, $C^{\rand}_{\res}$  and $C^{\quant}_{\res}$ which depends on $r$, $\rho$ and $D$, such that for any $\ve>0$ the $\ve$-complexity of median approximation problem satisfies
\[
\comp_{\ve,\abs}^{\det}(\tfrr)\leq C^{\det}_{\abs} \left(\frac{1}{\ve}\right)^{1/(r+\rho)}\, \log\frac{1}{\ve},
\]  
\[
\comp_{\ve,\abs}^{\rand}(\tfrr)\leq C^{\rand}_{\abs} \left(\frac{1}{\ve}\right)^{1/(r+\rho+1/2)}\, \log^2\frac{1}{\ve},
\]
\[  
\comp_{\ve,\abs}^{\quant}(\tfrr)\leq C^{\quant}_{\abs} \left(\frac{1}{\ve}\right)^{1/(r+\rho+1)}\, \log\frac{1}{\ve}\log\log\frac{1}{\ve},
\]  
\[
\comp_{\ve,\res}^{\det}(\frr)\leq
 C^{\det}_{\res} \left(\frac{1}{\ve}\right)^{1/(r+\rho)}\, \log\frac{1}{\ve},
\] 
\[
\comp_{\ve,\res}^{\rand}(\frr)\leq
 C^{\rand}_{\res} \left(\frac{1}{\ve}\right)^{1/(r+\rho+1/2)}\, \log^2\frac{1}{\ve}\log\log\frac{1}{\ve},
\] 
\[
\comp_{\ve,\res}^{\quant}(\frr)\leq
 C^{\quant}_{\res} \left(\frac{1}{\ve}\right)^{1/(r+\rho+1)}\, \log\frac{1}{\ve}\log\log\frac{1}{\ve},
\] 
\end{thm}  
\begin{pf}
To proof this theorem we construct  algorithms $\v^{\det}$, $\v^{\rand}$ and $\v^{\quant}$. These algorithms use perturbed  bisection method for solving nonlinear equation $G(x)=0$ for $G(x)=\intl_0^x f(t)dt$. For the positive parameter $\ve$ let $x_i$, $i=1,2,\ldots$ be the successive bisection points and $G_i$ be the approximation of $G(x_i)$ computed by optimal deterministic, randomized or quantum integration algorithm with precision $\ve$ with probability at least $1-\delta^{\rand/\quant}$ in the randomized and quantum setting, respectively, where $\delta^{\rand}=\ve^2/\lceil\log \ve^{-1}\rceil$ and $\delta^{\quant}=1/(4\lceil\log\ve^{-1}\rceil)$. That is
\begin{equation}\label{precision}
|G_i-G(x_i)|\leq \ve \textrm{ with probability at least } 1-\delta^{\rand/\quant}.
\end{equation}
The successive interval is chosen based on these approximate values $G_i$. We finish the bisection procedure when for some $i_0$, $|G_{i_0}|\leq \ve$ or $i_0=i_{max}:=\lceil\log \ve^{-1}\rceil$. The algorithms return the value $\hxi:=x_{i_0}$.\\
Notice that in all previous steps $i=1,2,\ldots,i_0-1$ we have $|G_i|>\ve$. So, due to (\ref{precision}) if $G_i\geq 0$ then $G(x_i)=G_i+G(x_i)-G(i)\geq G_i -\ve\geq 0$ (with probability at least $1-\delta^{\rand/\quant}$ in the randomized and quantum settings). Similarly, if $G_i<0$, then with the same probability we have $G(x_i)<0$. So, with probability at least $(1-\delta^{\rand/\quant})^{i_{max}}$ last bisection interval contains the median. So, if the algorithm finishes when $i_0=i_{max}$, then  with probability at least $(1-\delta^{\rand/\quant})^{i_{max}}$ we have $\xi-\hxi\leq \ve$. With the same probability, if the algorithm finishes when $|G_{i_0}|\leq\ve$, then
\[
|G(x_{i_0})|\leq |G(x_{i_0})-G_{i_0}|+|G_{i_0}|\leq 2\ve.
\]       
So, we have
\begin{equation}\label{precision_solution}
|G(\hxi)-G(\xi)|\leq 2\ve \textrm{ or } |\xi-\hxi|\leq\ve \textrm{ with probability at least } (1-\delta^{\rand/\quant})^{i_{max}}.
\end{equation}
We now state the error bounds. Consider first the residual criterion. Due to (\ref{precision_solution}), $|G(\hxi)-G(\xi)|\leq 2\ve$ or 
$
|G(\hxi)-G(\xi)|\leq \sup_{x\in[0,1]}|G'(x)|\,|\hxi-\xi|\leq D|\hxi-\xi|\leq D\ve$. Thus,
\[
e_{\res}(f,\v^{\det/\rand/\quant})\leq \max\{2,D\}\ve.
\]
This holds with probability at least $(1-\delta^{\rand})^{i_{max}}\geq 1-i_{max}\,\delta^{\rand}\geq 1-\ve^2$ in the randomized setting and with probability at least $(1-\delta^{\quant})^{i_{max}}\geq 1-i_{max}\,\delta^{\quant}\geq 3/4$ in the quantum settings. 

Hence we have
\begin{equation}\label{e_res_det}
\e^{\det}_{\res}(\frr,\v^{\det})\leq\max\{2,D\}\ve
\end{equation}
and
\begin{equation}\label{e_res_quant}
\e^{\quant}_{\res}(\frr,\v^{\quant})\leq\max\{2,D\}\ve.
\end{equation}
To get the error bound in the randomized setting, note that always $e_{\res}^{\omega}(f,\v^{\rand})\leq 1/2$. So we have in the randomized setting
\begin{multline*}
\left(\e^{\rand}_{\res}(\frr,\v)\right)^2=\sup_{f\in\frr}\intl_{\Omega}(\e^{\omega}_{\res}(f,\v^{\rand}))^2d\P(\omega)\\
=\sup_{f\in\frr}\Big(\intl_{\e^{\omega}_{\res}(f,\v^{\rand})\leq \max\{2,D\}\ve}(\e^{\omega}_{\res}(f,\v^{\rand}))^2d\P(\omega)+\\
\intl_{\e^{\omega}_{\res}(f,\v^{\rand})>\max\{2,D\} \ve}(\e^{\omega}_{\res}(f,\v^{\rand}))^2d\P(\omega)\Big)
\leq \sup_{f\in\frr}\left(\max\{2,D\}\ve^2+1/4\ve^2\right)
\end{multline*}
Thus,
\begin{equation}\label{e_res_rand}
e^{\rand}_{\res}(\frr,\v)/\leq(\max\{2,D\}+1/2)\,\ve.
\end{equation}
Let us pass to the absolute criterion. Let us remind that we have additional assumption, that $f(x)\geq \gamma$ for all $x\in[0,1]$. Note that $G(\hxi)-G(\xi)=f(c)(\hxi-\xi)$ for some $c\in[0,1]$. Hence, we have $|\hxi-\xi|\leq 1/\gamma |G(\hxi)-G(\xi)|$. So, due to (\ref{precision_solution}), $|\hxi-\xi|\leq \max\{1,2/\gamma\}\ve$ with probability at least $1-i_{max}\,\delta^{\rand/\quant}$ in the randomized and quantum settings. Hence, similarly as for the residual criterion, we have
\begin{equation}\label{e_abs_det}
\e^{\det}_{\abs}(\tfrr,\v^{\det})\leq \max\{1,2/\gamma\}\ve,
\end{equation}
\begin{equation}\label{e_abs_rand}
\e^{\rand}_{\abs}(\tfrr,\v^{\rand})\leq (\max\{1,2/\gamma\}+1/2)\ve
\end{equation}
and
\begin{equation}\label{e_abs_quant}
\e^{\quant}_{\abs}(\tfrr,\v^{\quant})\leq \max\{1,2/\gamma\}\ve.
\end{equation}
Let us state the cost bounds. In the deterministic setting, the cost of one bisection step for both criteria is of order 
$O\left(\left(1/\ve\right)^{1/(r+\rho)}\right)$.
The maximal number of steps is at most  $i_{max}=\lceil\log\ve^{-1}\rceil$. So, the total cost is
\begin{equation}\label{cost_det}
\cost^{\det}(\frr,\v^{\det})\asymp \cost^{\det}(\tfrr,\v^{\det})=O\left(\left(\frac{1}{\ve}\right)^{1/(r+\rho)}\log\ve^{-1}\right).
\end{equation}
In the randomized setting the cost of one bisection step is 
\begin{align*}
O\left(\left(\frac{1}{\ve}\right)^{1/(r+\rho+1/2)}\log(1/\delta^{\rand})\right)&=O\left(\left(\frac{1}{\ve}\right)^{1/(r+\rho+1/2)}\log\frac{\log\ve^{-1}}{\ve^2}\right)\\
&=O\left(\left(\frac{1}{\ve}\right)^{1/(r+\rho+1/2)}\log\ve^{-1}\right).
\end{align*}
So, the total cost in the randomized setting is
\begin{equation}\label{cost_rand}
\cost^{\rand}(\frr,\v^{\rand})\asymp \cost^{\rand}(\tfrr,\v^{\rand})=O\left(\left(\frac{1}{\ve}\right)^{1/(r+\rho)}\log^2\ve^{-1}\right).
\end{equation}
In the quantum setting, cost of one bisection step is
\[
O\left(\left(\frac{1}{\ve}\right)^{1/(r+\rho+1)}\log(1/\delta^{\quant})\right)=O\left(\left(\frac{1}{\ve}\right)^{1/(r+\rho+1)}\log\log\ve^{-1}\right).
\]
So, the total cost is of order
\begin{equation}\label{cost_quant}
\cost^{\quant}(\frr,\v)\asymp \cost^{\quant}(\tfrr,\v)=O\left(\left(\frac{1}{\ve}\right)^{1/(r+\rho+1)}\log\ve^{-1}\log\log\ve^{-1}\right).
\end{equation}
Comparing the error bounds (\ref{e_res_det}), (\ref{e_res_rand}), (\ref{e_res_quant}), (\ref{e_abs_det}), (\ref{e_abs_rand}), (\ref{e_abs_quant})  with the cost bounds (\ref{cost_det}), (\ref{cost_rand}), (\ref{cost_quant}) we get the desired complexity bounds. 
\end{pf}

\subsection{Lower complexity bounds}\label{sec:lower}
The following theorem presents the lower bounds on the complexity of problem of approximating the median in the deterministic, randomized and quantum settings.
\begin{thm}
There exist positive constants $c^{\det}_{\abs}$, $c^{\rand}_{\abs}$, $c^{\quant}_{\abs}$, $c^{\det}_{\res}$, $c^{\rand}_{\res}$, $c^{\quant}_{\res}$ and $\ve_0$, such that for any $\ve\in(0,\ve_0)$
\[
\comp_{\ve,\res}^{\det}(\frr(1))\geq c^{\det}_{\abs}\left(\frac{1}{\ve}\right)^{{1}/({r+\rho})},
\]
\[
\comp_{\ve,\res}^{\rand}(\frr(1))\geq c^{\rand}_{\abs}\left(\frac{1}{\ve}\right)^{{1}/({r+\rho+1/2})},
\]
\[
\comp_{\ve,\res}^{\quant}(\frr(1))\geq c^{\quant}_{\abs}\left(\frac{1}{\ve}\right)^{{1}/({r+\rho+1})},
\]
\[
\comp_{\ve,\abs}^{\det}(\tfrr(1))\geq c^{\det}_{\abs}\left(\frac{1}{\ve}\right)^{{1}/({r+\rho})},
\]
\[
\comp_{\ve,\abs}^{\rand}(\tfrr(1))\geq c^{\rand}_{\abs}\left(\frac{1}{\ve}\right)^{{1}/({r+\rho+1/2})},
\]
\[
\comp_{\ve,\abs}^{\quant}(\tfrr(1))\geq c^{\quant}_{\abs}\left(\frac{1}{\ve}\right)^{{1}/({r+\rho+1})},
\]
\end{thm}
\begin{pf}
We will proof the theorem above by reducing the problem of approximating the median of the distribution to the problem of approximating the mean of the discrete sequence.
Let $\ve_1>0$ be a parameter to be specified later on.  The class $\frr$ contains $n=\Theta\left(\ve_1^{-1/(r+\rho)}\right)$ functions $h_i$, $i=1,2,\ldots,n$, with disjoint supports in the interval $[0,1/4]$ and $g_i$, $i=1,2,\ldots,n$ with disjoint supports in the interval $[3/4,1]$ such that \[\intl_0^1h_i(x)dx=\intl_0^1g_i(x)dx=\ve_1^{1+1/(r+\rho)},\quad
\max_{x\in[0,1]}h_i(x)=\max_{x\in[0,1]}g_i(x)=c\,\ve_1,\quad i=1,2,\ldots n
\] 
for some constant $c$
 (see~\cite{Novak},~p.~35).
Let $x_1,x_2,\ldots,x_n$ be a sequence of real numbers in $[0,1]$. Then the function $\displaystyle f_{\ve_1}(x)=1+\sum_{i=1}^nx_ih_i(x)-\sum_{i=1}^nx_ig_i(x)$ belongs to the class $\tfrr(1)\subset\frr(1)$. \\
Note that for sufficiently small $\ve_1$ we have $2/3\leq f_{\ve_1}(x)\leq 2$ for all $x\in[0,1]$. Let $\xi$ be the median of $f_{\ve_1}$. Since $1/2=\int_0^{\xi}f_{\ve_1}(x)dx=\xi\,f_{\ve_1}(\eta)$ for some $\eta\in[0,1]$, thus $\xi=1/(2f_{\ve_1}(\eta))$ and $1/4\leq \xi\leq 3/4$. This yields that
\begin{align*}
1/2&=\intl_0^{\xi}f_{\ve_1}(x)dx=\intl_0^{\xi}\left(1+\sum_{i=1}^nx_ih_i(x)\right)dx=\intl_0^{1/4}\left(1+\sum_{i=1}^nx_ih_i(x)\right)dx+\intl_{1/4}^{\xi}1dx\\
&=1/4 + \ve_1^{1+1/(r+\rho)}\sum_{i=1}^n x_i+\xi-1/4.
\end{align*}
Thus
\[
\sum_{i=1}^nx_i=\frac{1/2-\xi}{\ve_1^{1+1/(r+\rho)}}.
\]
Suppose that algorithm $\v$ (deterministic, randomized or quantum) computes the median of the distribution with a density function $f$ with error at most $\ve$ and cost $N$, for any function $f\in\tfrr(1)$, in particular for $f=f_{\ve_1}$. Denote the result of the algorithm $\v$ for function $f_{\ve_1}$ by $\hat\xi$.  Note that 
\[
\left|\frac{1}{n}\sum_{i=1}^nx_i-\frac{1/2-\hat\xi}{n\ve_1^{1+1/(r+\rho)}}\right|=\frac{|\hat\xi-\xi|}{n\ve_1^{1+1/(r+\rho)}}=\frac{\left|\int_0^{\hat\xi}f_{\ve_1}(x)dx-1/2\right|}{|f_{\ve_1}(\eta)|n\ve_1^{1+1/(r+\rho)}}
\]
for some $\eta\in[0,1]$. Hence, algorithm $\v$ also computes the mean $\displaystyle\frac{1}{n}\sum_{i=1}^n x_i$ with error at most $\displaystyle\frac{C\ve}{n\ve_1^{1+1/(r+\rho)}}$, where $C=1$ for residual criterion and $C=1/\gamma$ when the absolute criterion is used. 

We now use lower complexity bounds for the problem of computing the mean of $n$ real numbers from $[0,1]$ presented in Section~\ref{sec:known}.\\
Consider first the deterministic setting. Recall that $n=\Theta\left(\varepsilon_1^{-1/(r+\rho)}\right)$. Let $G$ be a constant in the lower bound in the "$\Theta$" notation.
Take $\ve_1=4C\varepsilon/G$. From (\ref{lower_detst}) we have that the cost of algorithm $\v$ is bounded by
\begin{eqnarray*}
\displaystyle 
N&=&\Omega\left( n\left(1-\frac{2C\ve}{n\ve_1^{1+1/(r+\rho)}}\right)\right) =\Omega\left(n-\frac{2C\ve}{\ve_1^{1+1/(r+\rho)}}\right)\\
&=&\Omega\left(G\ve_1^{-1/(r+\rho)}-\frac{1}{2}G\ve_1\ve_1^{-1-1/(r+\rho)}\right)=\Omega\left(\ve^{-1/(r+\rho)}\right).
\end{eqnarray*}
\\ 
In the randomized setting, we take $\displaystyle\ve_1=\ve^{(r+\rho)/(r+\rho+1/2)}$. From (\ref{lower_rand}) we have
\[N=\Omega\left(\min\left\{n,\left(n\ve_1^{(r+\rho+1/2)/(r+\rho)}\ve^{-1}/C\right)^2\right\}\right) =\Omega\left(\ve^{-\frac{2}{2(r+\rho)+1}}\right)=\Omega\left(\ve^{-\frac{1}{r+\rho+1/2}}\right).\]  
\\
In the quantum setting, we take $\displaystyle\ve_1=\ve^{(r+\rho)/(r+\rho+1)}$. From (\ref{lower_quant}) we have
\[N=\Omega\left(\min\left\{n,n\ve_1^{(r+\rho+1)/(r+\rho)}\ve^{-1}/C\right\}\right)= \Omega\left(\ve^{-1/(r+\rho+1)}\right).\]  
The cost of any algorithm for computing the median must be at least $N$. This yields the desired lower bounds on the complexity in the deterministic, randomized and quantum settings.
\end{pf}
Note that, since the lower and the upper complexity bounds match up to the logarithmic factor, the algorithm presented in Section~\ref{sec:upper} is almost optimal.
\section{Computing  of quantiles}\label{sec:quantiles}
\subsection{Problem formulation}
In this section we consider the problem of approximation of the vector of quantiles of absolute continuous distribution. Let $f$ be density function. Suppose that for $k\in\N$ we are given a vector $\alpha=(\alpha_1,\alpha_2,\ldots,\alpha_k)\in[0,1]^k$. Our aim is to approximate the vector of quantiles $\xi=(\xi_1,\xi_2,\ldots,\xi_k)$, such that $\int_0^{\xi_1}f(x)dx=\alpha_1,\ \int_0^{\xi_2}f(x)dx=\alpha_2,\ldots,\int_0^{\xi_k}f(x)dx=\alpha_k$. We assume here that function $f$ is separated from zero and belongs to the H\:older class $\tfrr(1)$ defined in Section~\ref{sec:formulation}. We use the absolute error criterion, thus the local error of algorithm $\v$ approximating the vector of quantiles $\xi$ for the density function $f$ is defined by:
\begin{equation}
\e(f,\v)=\max_{i=1,\ldots,k}|\hat\xi_i-\xi_i|,
\end{equation}
where $\hat\xi=(\xi_1,\xi_2,\ldots,xi_k)$ is the approximation returned by the algorithm $\v$. The global error in the class $\tfrr$, the cost and the complexity are defined similarly as in Section~\ref{sec:formulation}. 
\subsection{Complexity bounds}
The following theorem presents the complexity bounds for the problem of approximation of the vector of the quantiles.

\begin{thm}
There exist positive constants $C^{\det}$, $C^{\rand}$, $C^{\quant}$, $c^{\det}$, $c^{\rand}$, $c^{\quant}$ and $\ve_0$ which depend on $r$, $\rho$, $D$ and $\gamma$, such that for any $k\in \N$, $\alpha\in[0,1]^k$ and $\ve>\ve_0$ the $\ve$-complexity of quantiles approximation problem satisfies
\[
c^{\det} \left(\frac{1}{\ve}\right)^{1/(r+\rho)}\leq\comp_{\ve,\abs}^{\det}(\tfrr)\leq C^{\det} \left(\frac{1}{\ve}\right)^{1/(r+\rho)},
\]  
\[
c^{\rand} \left(\frac{1}{\ve}\right)^{1/(r+\rho+1/2)}\leq\comp_{\ve,\abs}^{\rand}(\tfrr)\leq C^{\rand} \left(\frac{1}{\ve}\right)^{1/(r+\rho+1/2)},
\]
\[  
c^{\quant} \left(\frac{1}{\ve}\right)^{1/(r+\rho+1)}\leq\comp_{\ve,\abs}^{\quant}(\tfrr)\leq C^{\quant} \left(\frac{1}{\ve}\right)^{1/(r+\rho+1-\delta)},
\] 
\end{thm} 
\begin{pf}
Let $f\in\tfrr(1)$. Define the function $F:[0,1]\rightarrow[0,1]$ by $F(x)=\int_{0}^xf(x)dx$. Let $G=F^{-1}$ be the inverse of $F$. Then, $\xi=(\xi_1,\xi_2,\ldots,\xi_k)=(G(\alpha_1),G(\alpha_2),\ldots,G(\alpha_k))$. Note that $G$ is the solution of the initial-value problem
\[
G'(x)=\frac{1}{f(x)},\ \ x\in[0,1],\quad G(0)=0.
\]
Since $f\in\tfrr$, the function $1/f(x)$ is also H\"older continuous with the same parameters $r$ and $\rho$, but with different constant $\tilde D$ dependent on $D$, $r$ and $\gamma$. Let $l$ be the approximation of $G$ returned by the optimal algorithm solving initial-value problems in the H\"older class (deterministic, randomized or quantum). Then, the approximation of $\xi$ is defined by
\[
\tilde\xi=(\tilde\xi_1,\tilde\xi_2,\ldots,\tilde\xi_k)=(l(\alpha_1),l(\alpha_2),\ldots,l(\alpha_k)). 
\]
The upper bounds on the complexity of the problem of approximating of the vector of quantiles follow directly from the upper bounds on the problem of solving the initial-value problems presented in Section~\ref{sec:known}.

The lower bounds follow from the lower bounds on the complexity of approximating of the median presented in Section~\ref{sec:lower}. 
\end{pf}
\section{Remarks}
Note that the upper bounds for the problem of approximating the vector of quantiles improve slightly the bounds on the problem of approximating the median in the deterministic and randomized setting for the absolute error criterion. This bounds match the lower complexity bounds. All the other bounds are almost sharp.

\end{document}